\documentclass[10pt]{article}
\usepackage{amsfonts}
\usepackage{hyperref}
\usepackage{mathrsfs,enumerate,amssymb,amsmath,color,lineno}
\usepackage{indentfirst}
\usepackage{amsmath}
\usepackage{amssymb}
\usepackage[amsmath,thmmarks]{ntheorem}
\usepackage{cite}
\usepackage{graphicx}
\usepackage{tikz}

\newtheorem{definition}{\bf Definition}[section]
\newtheorem{lemma}{\bf Lemma}[section]
\newtheorem{theorem}{\bf Theorem}[section]
\newtheorem{remark}{\bf Remark}[section]
\newtheorem{corollary}{\bf Corollary}[section]
\newtheorem{example}{\bf Example}[section]

\newtheorem{proposition}{\bf Proposition}[section]

\begin{document}
\setcounter{page}{1}

\title{{\textbf{Monotone functions that generate conditionally cancellative triangular subnorms}}\thanks {Supported by
the National Natural Science Foundation of China (No.12471440)}}
\author{Meng Chen\footnote{\emph{E-mail address}: mathchen2019@163.com}, Xue-ping Wang\footnote{Corresponding author. xpwang1@hotmail.com; fax: +86-28-84761502},\\
\emph{School of Mathematical Sciences, Sichuan Normal University,}\\
\emph{Chengdu 610066, Sichuan, People's Republic of China}}

\newcommand{\pp}[2]{\frac{\partial #1}{\partial #2}}
\date{}
\maketitle
\begin{quote}
{\bf Abstract} Let a function $F: [0,1]^2\rightarrow [0,1]$ be given by
 $F(x,y)= f^{(-1)}(T(f(x), f(y)))$ where $f :[0,1]\rightarrow [0,1]$ is a monotone function, $f^{(-1)}$ is the pseudo-inverse of $f$ and $T$ is a triangular norm. This article characterizes the monotone function $f$ satisfying that the function $F$ is a conditionally cancellative triangular subnorm completely. It finally answers an open problem posed by Mesiarov\'{a}.

\end{quote}

\section{Introduction}

In 1826, Abel \cite{Abel} posed the problem: are there constructions involving only a one-place real 
function and the usual addition (or multiplication), which induce two-place real functions (i.e., operations on the real line) having nice algebraic properties, 
in particular, the associativity? Moreover, he constructed an associative function by using continuous strictly decreasing one-place functions and the usual addition. Following Abel's idea, Vicen\'{\i}k \cite{PV2005} gave the characterization of strictly monotone additive generators of associative functions. Mesiarov\'{a} \cite{AM2002} showed that if a continuous, Archimedean, proper t-subnorm is cancellative on $(0, 1]^2$ then it has a continuous additive generator. She  even characterized all continuous strictly monotone additive generators yielding continuous Archimedean proper t-subnorms, also mentioning an open problem \cite{AM2004}: Characterize all additive generators yielding continuous Archimedean proper t-subnorms? After that, she further studied continuous additive generators of continuous, conditionally cancellative t-subnorms, showing the conditions under which continuous, conditionally cancellative t-subnorms possess a continuous additive generator \cite{AM2016}. Recently, Zhang and Wang \cite{YM2024} supplied the characterization of one-place right continuous monotone functions that generate associative functions. Chen, Zhang and Wang \cite{2025chen} also presented the characterization of a class of one-place monotone functions which generate associative functions. 

On the one hand, obviously, the current literature as above have imposed some limiting conditions on the one-place functions, such as a continuous function, a strictly monotone function, etc. On the other hand, it should be pointed out that t-subnorms can be viewed as the most general semigroups yielding a t-norm when applying ordinal sums of semigroups \cite{Klement} and they serve as the basic building blocks for constructing t-norms through the ordinal sum method \cite{Jenei}. Therefore, understanding the structure and properties of t-subnorms is crucial for describing t-norms. These motivate us to consider an interesting problem: what is a characterization of one-place monotone functions that generate conditionally cancellative t-subnorms? This article will answer the above problem. Specific speaking, in this article we consider what is a characterization of all monotone functions $f: [0,1]\rightarrow [0,1]$ such that the function $F: [0,1]^2\rightarrow [0,1]$ given by \begin{equation}\label{eq:(5)}
F(x,y)= f^{(-1)}(T(f(x),f(y)))
\end{equation}
is a conditionally cancellative triangular subnorm, where $f^{(-1)}$ is the pseudo-inverse of $f$ and $T: [0,1]^2\rightarrow [0,1]$ is a t-norm.

The rest of this article is organized as follows. In Section 2 we recall some basic concepts and results. In Section 3 we entirely describe the characterization of a monotone function $f$ such that the function $F$ given by Eq.\eqref{eq:(5)} is a conditionally cancellative t-subnorm. A conclusion is drawn in Section 4, in which we answer the open problem posed by Mesiarov\'{a}.

\section{Preliminaries}

In this section, we recall some known basic concepts and results.
\begin{definition}[\cite{EP2000,CA2006}]
\emph{A \emph{triangular norm} (\emph{t-norm} for short) is a binary operator $T:[0, 1]^2\rightarrow [0, 1]$ such that for all $x, y, z\in[0, 1]$ the following conditions are satisfied:}

\emph{(T1)}  $T(x,y)=T(y,x)$,

\emph{(T2)}  $T(T(x,y),z)=T(x,T(y,z))$,

\emph{(T3)  $T(x,y)\leq T(x,z)$ whenever $y\leq z$,}

\emph{(T4)}  $T(x,1)=x$.

\emph{A binary operator $T:[0, 1]^2\rightarrow [0, 1]$ is called \emph{t-subnorm} if it satisfies (T1)-(T3), and $T(x,y)\leq \min\{x,y\}$ for all $x,y\in [0, 1]$.}
\end{definition}

Obviously, a t-norm is a t-subnorm. A binary operator $T:[0, 1]^2\rightarrow [0, 1]$ is called a \emph{proper t-subnorm} if it is a t-subnorm but not t-norm. 
\begin{definition}[\cite{EP2000,CA2006}]\label{de2.2}
	\emph{A \emph{triangular conorm} (\emph{t-conorm} for short) is a binary operator $S:[0, 1]^2\rightarrow [0, 1]$ such that for all $x, y, z\in[0, 1]$ the following conditions are satisfied:}
	
	\emph{$(S1)$ $S(x,y)=S(y,x)$,}
	
	\emph{$(S2)$  $S(S(x,y),z)=S(x,S(y,z))$,}
	
	\emph{$(S3)$  $S(x,y)\leq S(x,z)$ whenever $y\leq z$,}
	
	\emph{$(S4)$ $S(x,0)=x$.}
	
	\emph{A binary operator $S:[0, 1]^2\rightarrow [0, 1]$ is called a \emph{t-supconorm} if it satisfies $(S1), (S2), (S3)$, and $S(x,y)\geq \max\{x,y\}$ for all $x,y\in [0, 1]$.}
\end{definition}

\begin{definition}[\cite{EP2000,AM2016}]
\label{definition:2.3}
 \emph{For an arbitrary t-subnorm $T$ we consider the following properties:
 \begin{enumerate}
\renewcommand{\labelenumi}{(\roman{enumi})}
\item A t-subnorm $T$ is said to be \emph{strictly monotone} if $T(x,y)<T(x,z)$ whenever $x>0$ and $y<z$ for all $x,y,z\in[0,1]$.
\item The t-subnorm satisfies the \emph{cancellation law} if $T(x,y)=T(x,z)$ for all $x,y,z\in[0,1]$ implies $x=0$ or $y=z$;
\item The t-subnorm satisfies the \emph{conditional cancellation law} if for all $x,y,z\in[0,1]$, $T(x,y)=T(x,z)>0$ then $y=z$;
\item A t-subnorm $T$ is said to be \emph{Archimedean} if for every $(x,y)\in(0,1)^2$ there exists an $n\in\mathbb{N}$ such that $x_{T}^{(n)}<y$ where $x_{T}^{(n)}=T(x, x, \cdots, x)$ in which $x$ is presented $n$-times;
  \item A t-subnorm $T$ is said to be \emph{strict} if it is continuous and strictly monotone.
 \end{enumerate}}
\end{definition}

Denote
$$\mathcal{A}^*=\{M \mid \mbox{there is a monotone function } f:[0,1]\rightarrow[0,1]\mbox{ such that }\mbox{Ran}(f)=M\}$$
where $\mbox{Ran}(f)=\{f(x)|x\in [0,1]\}$.

Let $f:[0,1]\rightarrow[0,1]$ be a monotone function. Denote $f(a^-)=\lim_{x \rightarrow a^{-}}f(x)$ for all $a \in (0, 1]$ and, $f(a^+)=\lim_{x \rightarrow a^{+}}f(x)$ for all $a \in [0, 1)$. Define $f(0^{-}) = 0$ and $f(1^{+}) = 1$, respectively, whenever $f$ is a non-decreasing function, and $f(0^{-}) = 1$ and $f(1^{+}) = 0$, respectively, whenever $f$ is a non-increasing function. Denoted by $A\setminus B=\{x\in A\mid x\notin B\}$ for two sets $A$ and $B$.
\begin{lemma}\label{le2.1}
Let $M\in \mathcal{A}^*$ with $M \neq [0,1]$. Then there are a uniquely determined non-empty countable system $\mathcal{S}=\{[b_{k}, d_{k}] \subseteq [0, 1]\mid k\in K\}$ of closed intervals of a positive length which satisfy that for all $[b_{k}, d_{k}], [b_{l}, d_{l}]\in \mathcal{S}, [b_{k}, d_{k}]\cap[b_{l}, d_{l}]=\emptyset $ or $ [b_{k}, d_{k}]\cap[b_{l}, d_{l}]=\{d_{k}\} $ when $d_{k}\leq b_{l}$, and a uniquely determined non-empty countable set $C=\{c_{k} \in [0, 1] \mid k \in \overline{K}\}$ such that
$[b_{k}, d_{k}]\cap C\in\{\{b_{k}\}, \{d_{k}\}, \{b_{k}, d_{k}\}\}$ for all $k \in K$ and
\begin{equation*}
M= \{c_k\in [0,1]\mid k\in \overline{K}\}\cup \left([0,1]\setminus \left(\bigcup_{k\in K}[b_k, d_k] \right)\right)
\end{equation*}
where $|K|\leq|\overline{K}|$.
\end{lemma}
\begin{proof}
	 The proof is completely analogous to the proof of Lemma 3.1 of \cite{2025chen}.
\end{proof}

\begin{definition}\label{def2.1}
\emph{Let $M\in \mathcal{A}^*$. A pair $(\mathcal{S},C)$ is said to be \emph{associated with} $M\neq [0,1]$ if $\mathcal{S}=\{[b_{k}, d_{k}] \subseteq [0, 1]\mid k\in K \}$ is a non-empty system of closed intervals of a positive length which satisfy that for all $[b_{k}, d_{k}],[b_{l}, d_{l}]\in \mathcal{S}, [b_{k}, d_{k}]\cap[b_{l}, d_{l}]=\emptyset $ or $[b_{k}, d_{k}]\cap[b_{l}, d_{l}]=\{d_{k}\}$ when $d_{k}\leq b_{l}$, and $C=\{c_{k}\in [0, 1]\mid k\in \overline{K}\}$ is a non-empty countable set such that
$[b_{k}, d_{k}]\cap C\in\{\{b_{k}\}, \{d_{k}\}, \{b_{k}, d_{k}\}\}$ for all $k\in K$, and
\begin{equation*}
M= \{c_k\in [0,1]\mid k\in \overline{K}\}\cup \left([0, 1]\setminus \left(\bigcup_{k\in K}[b_k, d_k] \right)\right).
\end{equation*}
A pair $(\mathcal{S},C)$ is said to be \emph{associated with} $M=[0,1]$ if $\mathcal{S}=\{[1,1]\}$ and $C=\{1\}$.}
\end{definition}

\begin{definition}[\cite{EP2000,PV1999}]\label{def2.2}
\emph{Let $m,n,s, t\in [-\infty, \infty]$ with $m<n, s<t$ and $f:[m,n]\rightarrow[s,t]$ be an non-decreasing (resp. non-increasing) function. Then the function $f^{(-1)}:[s,t]\rightarrow[m,n]$ defined by
\begin{equation*}
f^{(-1)}(y)=\sup\{x\in [m,n]\mid f(x)< y\}\,(\mbox{resp. }f^{(-1)}(y)=\sup\{x\in [m,n]\mid f(x)> y\})
\end{equation*}
is called a \emph{pseudo-inverse} of the non-decreasing (resp. non-increasing) function $f$.}
\end{definition}

 If $f$ is a non-decreasing (resp. non-increasing) function, then we always have $f^{(-1)}\circ f\leq \mbox{id}_{[m,n]}$. Moreover, if $f$ is either right-continuous or strictly monotone then $f\circ f^{(-1)}\circ f=f$. In particular, let $M\in \mathcal{A}^*$ with $\mbox{Ran}(f)=M$ and $(\mathcal{S},C)=(\{[b_k, d_k]\mid k\in K\}, \{c_k \mid k\in K\})$ be associated with $M$. Then $f\circ f^{(-1)}(x)=x$ for all $x\in[0,1]$ if and only if $x\in M$. If $x\notin M$ then there is a $k\in K$ such that $x\in [b_k, d_k]\setminus\{c_{k}\}$ and $f\circ f^{(-1)}(x)=a_{k}$ where $a_{k}=b_{k}$ or $a_k=d_{k}$. In particular, if $f$ is a strictly monotone function, then $x\notin M$ implies that there is a $k\in K$ such that $x\in [b_k, d_k]\setminus\{c_{k}\}$ and $f\circ f^{(-1)}(x)=c_{k}$ where $\{c_{k}\}=M\cap [b_k, d_k]$

Then we have the following theorems.
\begin{theorem}[\cite{EP2000}]\label{th2.1}
For a function $T: [0, 1]^{2} \rightarrow [0,1]$ the following are equivalent:
\renewcommand{\labelenumi}{(\roman{enumi})}
\begin{enumerate}
\item $T$ is a strict t-norm.
\item  There exists a continuous,
strictly decreasing bijection $f : [0, 1] \rightarrow [0, \infty]$ such that
$T(x ,y) = f^{(-1)}(f(x) + f(y))$ for all $(x,y)\in [0, 1]^{2}$.
\end{enumerate}
\end{theorem}

 \begin{theorem}[\cite{AM2016}]\label{th2.2}
For a function $T: [0, 1]^{2} \rightarrow [0,1]$ the following are equivalent:
\renewcommand{\labelenumi}{(\roman{enumi})}
\begin{enumerate}
\item $T$ is a continuous, cancellative t-subnorm.
\item  There exists a continuous,
strictly decreasing function $f : [0, 1] \rightarrow [0, \infty]$ with $f(0)=\infty$ such that
$T(x ,y) = f^{(-1)}(f(x) + f(y))$ for all $(x,y)\in [0, 1]^{2}$.
\end{enumerate}
\end{theorem}
\begin{lemma}[\cite{AM2004}]\label{lemma2.1}
		A continuous t-subnorm $T$ is proper if and only if $T(1, 1)<1$.
\end{lemma}

From Corollary 3.1 (iii) of \cite{Wang2005} we have the following result.
\begin{proposition}\label{prop2.1}
	Let a function $F$ be given by Eq.~\eqref{eq:(5)} in which $f:[0,1]\rightarrow [0,1]$ is a strictly increasing function and $T$ is a strict t-norm.
	Then $\emph{Ran}(f)\cap[T(f(x^-),f(y^-)),T(f(x^+),f(y^+))]$ is at most a one-element set for all $x,y\in(0,1]$ if and only if $F$ is continuous.
\end{proposition}
\section{Generating functions of conditionally cancellative t-subnorms}
 In this section we characterize monotone functions $f:[0,1]\rightarrow [0,1]$ such that  $F: [0,1]^2\rightarrow [0,1]$ given by Eq.~\eqref{eq:(5)} is a conditionally cancellative t-subnorm completely.

For a monotone function $f:[0,1]\rightarrow [0,1]$, if the function $F: [0,1]^2\rightarrow [0,1]$ given by Eq.\eqref{eq:(5)} is a t-subnorm, then it is easy to check that $f$ is either a non-increasing function that satisfies $T(f(1),f(1))=f(1)$ with $f(0^+)=f(1)$ or a non-decreasing function. Therefore, in what follows, we always suppose that $f$ is a non-decreasing function and define $$Q=\{\omega\in[0,1]\mid \mbox{there exist two elements } x,y \in [0,1] \mbox{ with }x\neq y\mbox{ such that }f(x)=f(y)=\omega \}.$$
 Then one can verify the following proposition.
 \begin{proposition}\label{prop4.01}
Let a function $F$ be given by Eq.\eqref{eq:(5)} in which $f:[0,1]\rightarrow [0,1]$ is a non-decreasing function.  
\renewcommand{\labelenumi}{(\roman{enumi})}
\begin{enumerate}
\item If $f(1)\in Q$, then $F$ is a conditionally cancellative t-subnorm if and only if $F(x,y)=0$ for all $(x,y)\in[0,1]^{2}$.
\item If $f(1^{-})\in Q $, $f(1)\notin Q $, then $F$ is a conditionally cancellative t-subnorm if and only if 
\begin{equation*}
F(x,y)=\begin{cases}
a & \hbox{if }\ x=y=1,\\
0 & \hbox{otherwise}
\end{cases}
\end{equation*}
where $a\in [0,1]$.
\end{enumerate}
\end{proposition}

Therefore, in the rest of this section, we always suppose that $f:[0,1]\rightarrow [0,1]$ is a non-decreasing function with $f(1^-), f(1)\notin Q$.

Let $\upsilon=\sup Q$. Then $\upsilon=\max Q$ when $Q\neq \emptyset$, $\upsilon=0$ when $Q=\emptyset$. Hence there exists an element $\tau \in [0,1]$ such that $f(\tau^{-})=\upsilon$ and $f(\tau+\varepsilon)> \upsilon$ for an arbitrary $\varepsilon \in (0,1-\tau)$. Moreover, $f\mid_{[\tau,1]}$ is strictly increasing. Define $$\mathcal{N}=\{M \mid \mbox{there is a non-decreasing function } f:[0,1]\rightarrow[0,1]\mbox{ with } f(1^-),f(1)\notin Q \mbox{ such that }\mbox{Ran}(f)=M\}.$$ 
 Denoted by $T(A,B)=\{T(x,y)\mid x\in A, y\in B\}$ and $T(A,\emptyset)=\emptyset$. Then we have the following lemma.
\begin{lemma}\label{lem4.1}
Let $T:[0,1]^{2}\rightarrow [0,1]$ be a strictly monotone t-norm, $(\mathcal{S},C)$ be associated with $M \in \mathcal{N}$ and $\emph{Ran}(f)=M$. If $T(Q,M)\subseteq [0,f(0^{+})]$, then the following are equivalent:
\renewcommand{\labelenumi}{(\roman{enumi})}
\begin{enumerate}
\item $T(M\setminus C, M)\subseteq (M\setminus C) \cup [0,f(0^{+})]$.
\item $T(M\setminus C, M)\subseteq M \cup [0,f(0^{+})]$.
\end{enumerate}
\end{lemma}
\begin{proof}
 It is trivial that (i) implies (ii). We prove that (ii) implies (i). Let $T(M\setminus C,M)\subseteq M\cup [0,f(0^{+})]$. Then for all $x,y\in [0,1]$ with $f(x)\in M\setminus C$ we have $T(f(x),f(y))\in M \cup [0,f(0^{+})]$. If $T(f(x),f(y))\in [0,f(0^{+})]$ then clearly (i) is hold. If $T(f(x),f(y))\in M \setminus [0,f(0^{+})]$, then $x,y\in [\tau,1]$ since $T(Q,M)\subseteq [0,f(0^{+})]$. Suppose that $T(M\setminus C,M)\nsubseteq (M\setminus C) \cup [0,f(0^{+})]$.
Then there exist elements $x \in M\setminus C, y \in M$ and $k \in K$ such that $T(x,y) \in [b_{k}, d_{k}]$ with $b_{k} > 0$. Since $x \in (M\setminus C) \cap (0,1]$, there exists an $\varepsilon >0$ such that $a \in M\setminus C$ for any $a\in[x-\varepsilon, x]$. Thus there exists an $l\in K$ such that $T(a,y) \in [b_{l}, d_{l}]$ with $b_{l} > 0$ and, $T(a,y)\in M$. Moreover, from Definition \ref{def2.1} we have $M \cap [b_{l}, d_{l}]=\{c_{l}\}$, thus $T(a,y)= c_{l}$, in which $c_{l}$ has different values when $a$ changes since $T$ is a strictly monotone t-norm. Therefore, $\{[b_{l}, d_{l}]|l\in K\mbox{ and }c_{l}\in [b_{l}, d_{l}]\}$ is an infinitely uncountable set since $[x-\varepsilon, x]$ is infinitely uncountable, contrary to the fact that $K$ is a countable set.
\end{proof}

\begin{proposition}\label{prop4.1}
Let $T:[0,1]^{2}\rightarrow [0,1]$ be a strict t-norm, $(\mathcal{S},C)$ be associated with $M \in \mathcal{N}$ and $\emph{Ran}(f)=M$. Then the function $F$ given by Eq.~\eqref{eq:(5)} is conditionally cancellative if and only if the following hold:
\renewcommand{\labelenumi}{(\roman{enumi})}
\begin{enumerate}
\item $T(M\setminus C, M)\subseteq  M \cup [0,f(0^{+})]$ \mbox{and}
\item $T(Q,M)\subseteq  [0,f(0^{+})]$.
\end{enumerate}
\end{proposition}
\begin{proof}
We shall prove that $F$ is not
conditionally cancellative if and only if either $T(M\setminus C, M)\nsubseteq  M \cup [0,f(0^{+})]$ or $T(Q,M)\nsubseteq [0,f(0)]$.

Suppose that $F$ is not conditionally cancellative, i.e., there exist elements $x_{1}, x_{2}, y \in [0,1]$ with $x_{1} < x_{2}$ and $f(y)>0$ such that
$F(x_{1}, y) =F(x_{2},y) > 0$, i.e., $f^{(-1)}(T(f(x_{1}),f(y)))=f^{(-1)}(T(f(x_{2}),f(y)))>0$.
 In the following we distinguish two cases.
 \renewcommand{\labelenumi}{\rm $\centerdot$}
\begin{enumerate}
\item If $T(Q,M)\subseteq  [0,f(0^{+})]$, then $f(x_{1}), f(x_{2}), f(y)\notin Q$ since $F(x_{1}, y) =F(x_{2},y) > 0$. Thus, $T(f(x_{1}),f(y))< T(f(x_{2}),f(y))$ since $T$ is a strict t-norm and $f(y)>0$. Hence 
\begin{equation}\label{1002}
	[T(f(x_{1}),f(y)),T(f(x_{2}),f(y))] \cap (M\setminus C) = \emptyset.
\end{equation}
 On the other hand, $(M\setminus C)\cap(f(x_{1}), f(x_{2}))\neq \emptyset$ since $f(x_{1}), f(x_{2})\notin Q$, which implies that there is an $x \in [0,1]$ such that $f(x) \in (M\setminus C)\cap(f(x_{1}), f(x_{2}))$. This follows that $$T(f(x),f(y))\in [T(f(x_{1}),f(y)),T(f(x_{2}),f(y))]$$ which together with Eq. \eqref{1002} implies $T(f(x),f(y))\notin (M\setminus C)\cup [0,f(0^{+})]$. Thus, from Lemma \ref{lem4.1} we have $T(M\setminus C, M)\nsubseteq M\cup [0,f(0^{+})]$.

\item If $T(M\setminus C, M)\subseteq  M \cup [0,f(0^{+})]$ then, obviously, $T(Q,M)\nsubseteq [0,f(0^{+})]$.
\end{enumerate}

Conversely, suppose that $T(M\setminus C, M)\nsubseteq  M \cup [0,f(0^{+})]$, i.e., there exist elements $x,y \in [0,1]$ and $k\in K$ such that $f(x)\in M\setminus C$, $f(y) \in M$ and $T(f(x),f(y)) \in [b_{k}, d_{k}]\setminus\{c_{k}\}$ where $\{c_{k}\}=M\cap[b_{k}, d_{k}]$ and $b_{k}> f(0^{+})$. Because of $f(x) \in (M\setminus C) \cap [0,1]$, there exist $f(x_{1}), f(x_{2}) \in M$ with $f(x_{1}) \neq f(x_{2})$ such
that $$T(f(x_{1}), f(y)), T(f(x_{2}), f(y)) \in [b_{k}, d_{k}]\setminus\{c_{k}\}.$$ Then $F(x_{1}, y)=F(x_{2}, y)>0$. Thus $F$ is not conditionally cancellative. Now, suppose that $T(Q,M)\nsubseteq [0,f(0^{+})]$. Then there exist elements $x_{1}, x_{2},y\in [0,1]$ with $f(x_{1})=f(x_{2})$ such that $$T(f(x_{1}),f(y))= T(f(x_{2}),f(y))>f(0^{+}).$$ Therefore, $F(x_{1}, y)=F(x_{2}, y)>0$, i.e., $F$ is not conditionally cancellative.
\end{proof}

Notice that in Proposition \ref{prop4.1} the condition that $T$ is strict cannot be deleted generally.
\begin{example}\label{exp03.1}\emph{Let $T:[0,1]^{2}\rightarrow [0,1]$ be a t-norm, $(\mathcal{S},C)$ be associated with $M \in \mathcal{N}$ and $\mbox{Ran}(f)=M$.}
\renewcommand{\labelenumi}{(\roman{enumi})}
\emph{\begin{enumerate}
\item Let the function $f:[0,1]\rightarrow [0,1]$ be defined by
\begin{equation*}
 f(x)=\begin{cases}
\frac{1}{4}x+\frac{1}{4} & \hbox{if }\ x\in[0,1),\\
1 & \hbox{if }\ x=1
\end{cases}
\end{equation*}
and $T:[0,1]^{2}\rightarrow [0,1]$ be defined by
$T(x,y)=\min\{x,y\}$ for all $(x,y)\in[0,1]^{2}$.
Then by Eq.~\eqref{eq:(5)}
$$F(x,y)=\min\{x,y\} \mbox{ for all } (x,y)\in[0,1]^{2}.$$
Obviously, $T$ is a continuous but not strictly monotone t-norm, $T(M\setminus C, M)\subseteq  M \cup [0,f(0^{+})]$ and $T(Q,M)\subseteq  [0,f(0^{+})]$. However, $F$ is not conditionally cancellative.
\item Let the function $f:[0,1]\rightarrow [0,1]$ be defined by
\begin{equation*}
 f(x)=\begin{cases}
\frac{1}{2}x & \hbox{if }\ x\in[0,\frac{1}{4}),\\
\frac{13}{12}x-\frac{1}{12} & \hbox{if }\ x\in[\frac{1}{4},1]
\end{cases}
\end{equation*}
and $T:[0,1]^{2}\rightarrow [0,1]$ be defined by
\begin{equation*}
T(x,y)=\begin{cases}
\frac{xy}{2} & \hbox{if }\ (x,y)\in[0,\frac{1}{2}]^{2},\\
xy & \hbox{otherwise}.
\end{cases}
\end{equation*}
Then by Eq.~\eqref{eq:(5)}
\begin{equation*}
F(x,y)=\begin{cases}
\frac{xy}{4} & \hbox{if }\ (x,y)\in[0,\frac{1}{4})^{2},\\
\frac{13xy-x}{12} & \hbox{if }\ (x,y)\in[0,\frac{1}{4})\times [\frac{1}{4},1],\\
\frac{13xy-y}{12} & \hbox{if }\ (x,y)\in [\frac{1}{4},1]\times[0,\frac{1}{4}),\\
\frac{169xy-13x-13y+1}{144} & \hbox{if }\ (x,y)\in [\frac{1}{4},\frac{7}{13})^{2},\\
\frac{13xy-x-y+1}{12} & \hbox{otherwise}.
\end{cases}
\end{equation*}
It is easy to see that $T$ is a strictly monotone but non-continuous t-norm and $F$ satisfies the conditional cancellation law. One check that  $M=[0,\frac{1}{8})\cup[\frac{3}{16},1]$ and $\frac{1}{2}\in M\setminus C$. Moreover, $T(\frac{1}{2},\frac{1}{2})=\frac{1}{8}$. Thus 
$T(M\setminus C, M)\nsubseteq  M \cup [0,f(0^{+})]$.
\end{enumerate}}
\end{example}
Thus, in the rest of this article, we always suppose that $T:[0,1]^{2}\rightarrow [0,1]$ is a strict t-norm and define $$\mathcal{A}=\{M \mid \mbox{ there is a strictly increasing function } f:[0,1]\rightarrow[0,1] \mbox{ such that }\mbox{Ran}(f)=M\}.$$
In particular, if $f:[0,1]\rightarrow[0,1]$ is a strictly increasing function, then $Q=\emptyset$, thus $T(Q,M)=\emptyset\subseteq  [0,f(0^{+})]$. Therefore, from Proposition \ref{prop4.1}, we have the following corollary.
\begin{corollary}\label{cor3.01}
Let $T:[0,1]^{2}\rightarrow [0,1]$ be a strict t-norm, $(\mathcal{S},C)$ be associated with $M \in \mathcal{A}$ and $\emph{Ran}(f)=M$. Then the following are equivalent:
\renewcommand{\labelenumi}{(\roman{enumi})}
\begin{enumerate}
\item The function $F$ given by Eq.~\eqref{eq:(5)} is conditionally cancellative.
\item $T(M\setminus C, M)\subseteq M \cup [0,f(0^{+})]$.
\end{enumerate}
\end{corollary}

\begin{lemma}\label{lem04.1}
Let $T: [0,1]^2\rightarrow [0,1]$ be a strict t-norm and $f:[0,1]\rightarrow [0,1]$ be a non-decreasing function. If the function $F: [0,1]^2\rightarrow [0,1]$ given by Eq.\eqref{eq:(5)} satisfies the cancellation law, then $f$ is a strictly increasing function with $f(0)=f(0^{+})=0$.
\end{lemma}
\begin{proof}Assume that $f$ is not strictly increasing function. Then there exist two elements $x,y\in(0,1)$ with $x\neq y$ such that $f(x)=f(y)$. Thus, for any $z\in(0,1]$, $T(f(x),f(z))=T(f(y),f(z))$. By Eq. \eqref{eq:(5)}, we obtain 
$F(x,z)=F(y,z)$, a contradiction. Since the function $F: [0,1]^2\rightarrow [0,1]$ given by Eq. \eqref{eq:(5)} satisfies the cancellation law and $F(x,0)=0$ for all $x\in [0,1]$, we have 
$T(f(x),f(y))\geq f(0^+)$ for all $x,y\in (0,1]$, which implies $T(f(x),f(0^+))=f(0^+)$ for all $x\in (0,1]$. This follows that $f(0^+)=0$ since $T$ is a strict t-norm.
\end{proof}

 Therefore, we have the following proposition.
\begin{proposition}\label{prop3.2}
Let $T:[0,1]^{2}\rightarrow [0,1]$ be a strict t-norm, $(\mathcal{S},C)$ be associated with $M \in \mathcal{N}$ and $\emph{Ran}(f)=M$. Then the following are equivalent:
\renewcommand{\labelenumi}{(\roman{enumi})}
\begin{enumerate}
\item The function $F$ given by Eq.\eqref{eq:(5)} satisfies the cancellation law.
\item $T(M\setminus C,M)\subseteq M$.
\end{enumerate}
\end{proposition}
\begin{proof}
$(i\Rightarrow ii)$ If the function $F$ given by Eq.~\eqref{eq:(5)} satisfies the cancellation law then, from Lemma \ref{lem04.1} and Definition \ref{definition:2.3}, $f$ is a strictly increasing function with $f(0)=f(0^{+})=0$ and $F$ satisfies the conditional cancellation law. From Corollary \ref{cor3.01}, $T(M\setminus C, M)\subseteq M \cup [0,f(0^{+})]$. This follows that $T(M\setminus C,M)\subseteq M$ since $f$ is a strictly increasing function with $f(0)=f(0^{+})=0$.

$(ii\Rightarrow i)$ From Lemma \ref{lem04.1}, $f$ is a strictly increasing function with $f(0)=f(0^{+})=0$, then $T(Q,M)=\emptyset\subseteq  [0,f(0^{+})]$ and $[0,f(0^{+})]=\{0\}$. Thus, from Lemma \ref{lem4.1} $T(M\setminus C, M)\subseteq (M\setminus C) \cup \{0\}$ if and only if $T(M\setminus C, M)\subseteq M $. Suppose that $F$ does not satisfy the cancellation law, i.e., there exist elements $x_{1}, x_{2}, y \in [0,1]$ with $y > 0$ and $x_{1} < x_{2}$ such that
$F(x_{1}, y) =F(x_{2},y)$, i.e., $$f^{(-1)}(T(f(x_{1}),f(y)))=f^{(-1)}(T(f(x_{2}),f(y))).$$ Then $f(y)>0$ since $f$ is strictly increasing and, $T(f(x_{1}),f(y))< T(f(x_{2}),f(y))$ since $T$ is a strict t-norm. Thus 
\begin{equation}\label{eq001}
	[T(f(x_{1}),f(y)),T(f(x_{2}),f(y))] \cap (M\setminus C) = \emptyset.
\end{equation} 
On the other hand, we have $(M\setminus C)\cap(f(x_{1}), f(x_{2}))\neq \emptyset$ since $f$ is strictly increasing. Thus there exists an $x \in [0,1]$ such that $f(x) \in (M\setminus C)\cap(f(x_{1}), f(x_{2}))$. This follows that $$T(f(x),f(y))\in [T(f(x_{1}),f(y)),T(f(x_{2}),f(y))],$$ which together with Eq.~\eqref{eq001} implies that $T(f(x),f(y))\notin M\setminus C$. Therefore, $T(M\setminus C, M)\nsubseteq (M\setminus C) \cup \{0\}$, i.e., $T(M\setminus C,M)\nsubseteq M$, a contradiction.
\end{proof}

 \begin{lemma}\label{lem3.2}
Let $T:[0,1]^{2}\rightarrow [0,1]$ be a strict t-norm, $(\mathcal{S},C)$ be associated with $M \in \mathcal{A}$ and $\emph{Ran}(f)=M$. Then $f^{(-1)}(x)=f^{(-1)}(y)$ for all $x,y\in[0,1]$ if and only if $(\min\{x,y\},\max\{x,y\})\cap (M\setminus C)=\emptyset$.
\end{lemma}
\begin{proof}
If $f^{(-1)}(x)=f^{(-1)}(y)$ for all $x,y\in[0,1]$, then $f(f^{(-1)}(x))= f(f^{(-1)}(y))$. Suppose $$(\min\{x,y\},\max\{x,y\})\cap (M\setminus C)\neq\emptyset.$$ Then $x\neq y$, say $x<y$, thus there is an $a \in(x,y)\cap (M\setminus C)$. By Definition \ref{def2.2}, $f(f^{(-1)}(x))\leq x< f(f^{(-1)}(a))< f( f^{(-1)}(y))$, a contradiction.

Conversely, we know that $x\leq y$ or $y\leq x$ for all $x,y\in[0,1]$, say $x\leq y$. If $x=y$ then clearly $f^{(-1)}(x)=f^{(-1)}(y)$. If $x<y$, then from $(\min\{x,y\},\max\{x,y\})\cap (M\setminus C)=\emptyset$, there is a $k\in K$ such that $(x,y)\subseteq [b_k,d_k]$. Thus $f^{(-1)}(x)=f^{(-1)}(c_{k})=f^{(-1)}(y)$ where $\{c_{k}\} = M\cap[b_{k}, d_{k}]$.
\end{proof}

Let $(\mathcal{S},C)$ be associated with $M \in \mathcal{A}$. Define $$O(M) = \bigcup_{x,y\in M}(\min\{x,y\}, \max\{x,y\})$$ when $M\neq \emptyset$ and, $O(M)=\emptyset$ when $M=\emptyset$. Note that $(x,x)=\emptyset$ for any $x\in[0,1]$ and $T(c,O(M))=O(T(c,M))$ for any strict t-norm $T$. For each $k \in K$, define $$H_{k}=O(\{c_{k}\}\cup \{z\in[b_{k},d_{k}]\}\mid  \mbox{there exist two elements } x,y\in M \mbox{ such that } T(x,y)=z\}).$$ Then we have the following proposition.

\begin{proposition}\label{prop3.3}
 Let $T:[0,1]^{2}\rightarrow [0,1]$ be a strict t-norm, $(\mathcal{S},C)$ be associated with $M \in \mathcal{A}$ and $\emph{Ran}(f)=M$. If $T(\cup_{k\in K}H_{k}, M\setminus\{0\})\cap (M\setminus C)=\emptyset$, then the function $F$ given by Eq.~\eqref{eq:(5)} is a t-subnorm.
\end{proposition}
\begin{proof}
 It is enough to show that $F$ is associative. It is easy to see that $T(\cup_{k\in K}H_{k}, M\setminus\{0\})\cap (M\setminus C)= \emptyset$ implies 
 \begin{equation}\label{eq0002}
 	T(H_{k},f(z))\cap (M\setminus C)= \emptyset
 \end{equation} 
for any $k\in K$ and $z\in [0,1]$ with $f(z)\in(0,1] $. Next, we first prove that for any $x,y,z\in [0,1]$ with $f(x), f(y), f(z)\in(0,1]$, $$F(F(x,y),z)=f^{(-1)}(T((T(f(x),f(y))),f(z))).$$ We distinguish two cases as follows.

(i) If $T(f(x),f(y))\in M$, then $f(f^{(-1)}(T(f(x),f(y))))=T(f(x),f(y))$. It follows from Eq. \eqref{eq:(5)} that
\begin{align*}
 F(F(x,y),z) &= f^{(-1)}(T(f(f^{(-1)}(T(f(x),f(y))),f(z))))\\
 &=  f^{(-1)}(T((T(f(x),f(y))),f(z))).
\end{align*}

(ii) If $T(f(x),f(y))\notin M$, then there exists a $k \in K$ such that $T(f(x),f(y))\in [b_{k}, d_{k}]\setminus\{c_{k}\}$ where $\{c_{k}\} = M\cap[b_{k}, d_{k}]$. This together with Eq. \eqref{eq0002} implies that
$$T((\min\{c_{k}, T(f(x),f(y))\},\max\{c_{k}, T(f(x),f(y))\}), f(z))\cap (M\setminus C)= \emptyset$$
since $(\min\{c_{k}, T(f(x),f(y))\},\max\{c_{k}, T(f(x),f(y))\}) \in H_{k}$.
Thus by Lemma \ref{lem3.2},
$$f^{(-1)}(T(T(f(x),f(y))),f(z))= f^{(-1)}(T(c_{k},f(z))),$$
which deduces
$$ f^{(-1)}(T(f(f^{(-1)}T(f(x),f(y)),f(z))))=f^{(-1)}(T(c_{k},f(z)))=f^{(-1)}(T((T(f(x),f(y))),f(z)))$$ since $c_{k}=f(f^{(-1)}(T(f(x),f(y))))$.
Hence, $$F(F(x,y),z)=f^{(-1)}(T((T(f(x),f(y))),f(z))).$$

Analogously, we have $$F(x,F(y,z))=f^{(-1)}(T(f(x),T(f(y),f(z))))$$ for all $x,y,z\in [0,1]$ with $f(x), f(y), f(z)\in(0,1]$. Hence $F(F(x,y),z)=F(x,F(y,z))$ for all $x,y,z\in [0,1]$ with $f(x), f(y), f(z)\in(0,1]$ since $T$ is associative.
Moreover, if $0\in\{f(x), f(y), f(z)\}$, then clearly $F(F(x,y),z)=0=F(x,F(y,z))$. Therefore, $F(F(x,y),z)=F(x,F(y,z))$ for all $x,y,z\in [0,1]$, i.e., $F$ is associative.
\end{proof}

From Propositions \ref{prop3.2} and \ref{prop3.3}, the following  two corollaries are immediate.
\begin{corollary}\label{cor.3.2}
Let $T:[0,1]^{2}\rightarrow [0,1]$ be a strict t-norm, $(\mathcal{S},C)$ be associated with $M \in \mathcal{A}$ and $\emph{Ran}(f)=M$. If $T(M,M)\subseteq M\cup [0, f(0^{+}]$ and $T(M\setminus C,M)\subseteq M$, then the function $F$ given by Eq.~\eqref{eq:(5)} is a cancellative t-subnorm.
\end{corollary}

\begin{corollary}\label{cor.3.3}
Let $T:[0,1]^{2}\rightarrow [0,1]$ be a strict t-norm, $(\mathcal{S},C)$ be associated with $M \in \mathcal{A}$ and $\emph{Ran}(f)=M$. If $f$ is a strictly increasing continuous function with $f(0)=0$, then the function $F$ given by Eq.~\eqref{eq:(5)} is a cancellative t-subnorm.
\end{corollary}

Let $(\mathcal{S},C)$ be associated with $M \in \mathcal{N}$ and $\mbox{Ran}(f)=M$.  Define $$K_1=\{k\in K\mid b_k\geq \max Q\}.$$
\begin{proposition}\label{prop4.2}
Let $T:[0,1]^{2}\rightarrow [0,1]$ be a strict t-norm, $(\mathcal{S},C)$ be associated with $M \in \mathcal{N}$ and $\emph{Ran}(f)=M$. If the following statements hold:
\renewcommand{\labelenumi}{(\roman{enumi})}
\begin{enumerate}
\item $T(\cup_{k\in K_1}H_{k}, M\setminus\{0\})\cap (M\setminus C)= \emptyset$,
\item $T(M\setminus C, M)\subseteq M \cup [0,f(0^{+})]$ \mbox{ and}
\item $T(Q,M)\subseteq [0,f(0^{+})]$,
\end{enumerate}
 then the function $F$ given by Eq.~\eqref{eq:(5)} is a conditionally cancellative t-subnorm.
\end{proposition}
\begin{proof}
From Proposition \ref{prop4.1}, we just need to show that the function $F$ given by Eq.~\eqref{eq:(5)} is associative. Because $f\mid_{[\tau,1]}$ is strictly increasing, from Proposition \ref{prop3.3} we have $$F(F(x,y),z)=F(x,F(y,z))$$ for all $x,y,z\in(\tau,1]$. If $x\in[0,\tau]$, then $f(x)\leq \max Q$, which implies
$$F(F(x,y),z)=0=F(x,F(y,z))$$ for all $y,z\in[0,1]$ since $T(Q,M)\subseteq [0,f(0^{+})]$. If $y\in[0,\tau]$ or $z\in[0,\tau]$, then in complete analogy to $x\in[0,\tau]$ we have $F(F(x,y),z)=F(x,F(y,z))$ for all $x,z\in[0,1]$ or $F(F(x,y),z)=F(x,F(y,z))$ for all $x,y\in[0,1]$. Therefore, $F$ is associative.
\end{proof}

The following two corollaries are deduced from Propositions \ref{prop4.2}.
\begin{corollary}
Let $T:[0,1]^{2}\rightarrow [0,1]$ be a strict t-norm, $(\mathcal{S},C)$ be associated with $M \in \mathcal{N}$ and $\emph{Ran}(f)=M$. If the following statements hold:
\renewcommand{\labelenumi}{(\roman{enumi})}
\begin{enumerate}
\item $T(M, M)\subseteq M \cup [0,f(0^{+})]$ \mbox{ and}
\item $T(Q,M)\subseteq  [0,f(0^{+})]$,
\end{enumerate}
then the function $F$ given by Eq.\eqref{eq:(5)} is a conditionally cancellative t-subnorm.
\end{corollary}

\begin{corollary}
Let $T:[0,1]^{2}\rightarrow [0,1]$ be a strict t-norm, $(\mathcal{S},C)$ be associated with $M \in \mathcal{A}$ and $\emph{Ran}(f)=M$. If $T(M, M)\subseteq M \cup [0,f(0^{+})]$,
then the function $F$ given by Eq.~\eqref{eq:(5)} is a conditionally cancellative t-subnorm.
\end{corollary}

Generally, the converse of Proposition \ref{prop4.2} isn't true.
\begin{example}\label{exp3.01}
\emph{Let the function $f:[0,1]\rightarrow [0,1]$ be defined by
\begin{equation*}
 f(x)=\begin{cases}
\frac{1}{4}+ \frac{1}{4}x& \hbox{if }\ x\in[0,\frac{1}{4}),\\
\frac{5}{16} & \hbox{if }\ x\in[\frac{1}{4},\frac{1}{2}],\\
\frac{3}{8}+ \frac{1}{8}x& \hbox{if }\ x\in(\frac{1}{2},1),\\
\frac{3}{4} & \hbox{if }\ x=1
\end{cases}
\end{equation*}
and $T:[0,1]^{2}\rightarrow [0,1]$ be defined by
$T(x,y)=xy$ for all $(x,y)\in[0,1]^{2}$.
Then by Eq.~\eqref{eq:(5)}
\begin{equation*}
F(x,y)=\begin{cases}
1 & \hbox{if }\ xy=1,\\
0 & \hbox{otherwise }.
\end{cases}
\end{equation*}
Obviously, $f$ is a non-decreasing function with $\hbox{Ran}(f)=[\frac{1}{4},\frac{5}{16}]\cup (\frac{7}{16},\frac{1}{2})\cup \{\frac{3}{4}\}$ and $T$ is a strict t-norm. One can check that $F$ is a conditionally cancellative t-subnorm. However, $T(\cup_{k\in K_1}H_{k}, M\setminus\{0\})\cap (M\setminus C)= (\frac{1}{4},\frac{9}{32})$.}
\end{example}

Fortunately, we have the following proposition.
\begin{proposition}\label{th4.1}
Let $T:[0,1]^{2}\rightarrow [0,1]$ be a strict t-norm, $(\mathcal{S},C)$ be associated with $M \in \mathcal{N}$ and  $\emph{Ran}(f)=M$. If $f(1)=1$, then the function $F$ given by Eq.~\eqref{eq:(5)} is a conditionally cancellative t-subnorm if and only if the following statements hold:
\renewcommand{\labelenumi}{(\roman{enumi})}
\begin{enumerate}
\item $T(\cup_{k\in K_1}H_{k}, M\setminus\{0\})\cap (M\setminus C)= \emptyset$,
\item $T(M\setminus C, M)\subseteq M \cup [0,f(0^{+})]$ \mbox{and}
\item $T(Q,M)\subseteq [0,f(0^{+})]$.
\end{enumerate}
\end{proposition}
\begin{proof}
From Propositions \ref{prop4.1} and \ref{prop4.2}, we just need to show that if the function $F$ given by Eq.\eqref{eq:(5)} is a conditionally cancellative t-subnorm then $T(\cup_{k\in K_1}H_{k}, M\setminus\{0\})\cap (M\setminus C)= \emptyset$.

If $T(M , M) \subseteq M \cup [0,f(0^{+})]$ then either $\cup_{k\in K_1}H_{k}\subseteq[0,f(0^{+})]$ or $\cup_{k\in K_1}H_{k}=\emptyset$. If $\cup_{k\in K_1}H_{k}\subseteq[0,f(0^{+})]$ then $T(\cup_{k\in K_1}H_{k}, M\setminus\{0\}) \subseteq[0,f(0^{+})]$ since $T(x,y)\leq\min\{x,y\}$ for all $(x,y)\in[0,1]^{2}$. If $\cup_{k\in K_1}H_{k}=\emptyset$ then $T(\cup_{k\in K_1}H_{k}, M\setminus\{0\})=\emptyset$. Therefore, we always have $T(\cup_{k\in K_1}H_{k}, M\setminus\{0\})\cap (M\setminus C)= \emptyset$. Now, suppose that $T(M, M) \nsubseteq M \cup [0,f(0^{+})]$. Then there exist two elements $x,y \in [0,1]$ such that $T(f(x),f(y)) \notin M \cup [0,f(0^{+})]$, which means that there exists a $k \in K_1$ such that
\begin{equation}\label{eq0010}
	T(f(x),f(y))\in [b_{k}, d_{k}]\setminus\{c_{k}\}
\end{equation}
 where $\{c_{k}\}=[b_{k}, d_{k}] \cap M$. Moreover, we claim that $ f(1)>\max\{f(x),f(y)\}>0$ and $f(x),f(y)\notin Q$. Indeed, if $f(x)=f(y)=0$, then $0\in M$, thus $T(f(x),f(y))=T(0,0)=0\in M$, a contradiction. If $\max\{f(x),f(y)\}=f(1)$ then $T(f(x),f(y))=\min\{f(x),f(y)\}\in M$, a contradiction. If $f(x)\in Q$ or $f(y)\in Q$, then $T(f(x),f(y))\in [0,f(0^{+})]$, a contradiction.
Thus $\sup\{M\setminus\{1\}\} > \max\{f(x),f(y)\}$ since $f\mid_{(\tau,1]}$ is strictly increasing. This follows that there is a $z \in (\tau,1]$ such that 
\begin{equation}\label{eq0011}
	f(z)\in M\setminus C
\end{equation}
 with $f(z) > \max\{f(x),f(y)\}$.
On the other hand, by Proposition \ref{prop4.1} and Lemma \ref{lem4.1} we have
$$T(M\setminus C, M)\subseteq (M\setminus C)\cup [0,f(0^{+})].$$
Hence, from Eq. \eqref{eq0011} we have $T(f(y),f(z)) \in (M\setminus C)\cup [0,f(0^{+})]$. This follows $T(f(y),f(z))\in M\setminus C$ since $T(f(y),f(z))>T(f(y),f(x))>f(0)$, and so
$F(x,F(y,z))=f^{(-1)}T(f(x),T(f(y),f(z)))$. Meanwhile, Eqs. \eqref{eq:(5)} and \eqref{eq0010} yield $F(F(x,y),z)=f^{(-1)}(T(c_{k}, f(z)))$.
Thus, the associativity of $F$ and $T$ results in
$$f^{(-1)}(T(T(f(x), f(y)),f(z)))=f^{(-1)}(T(c_{k}, f(z))),$$
implying $$M\cap(\min\{T(T(f(x), f(y)),f(z))),T(c_{k}, f(z))\},\max\{T(T(f(x), f(y)),f(z))),T(c_{k}, f(z))\})=\emptyset.$$
If $T(c_{k}, f(z))> f(0^{+})$ then \ $T(c_{k}, f(z))\in M\setminus C$, implying $$T(T(f(x), f(y)),f(z)))=T(c_{k}, f(z)).$$ So that $T(f(x), f(y))=c_{k}$ since $T$ is a strictly monotone and $f(z) > f(0^{+})$, contrary to $T(f(x), f(y))\in [b_{k}, d_{k}]\setminus\{c_{k}\}$. Similarly, one can check that $T(T(f(x), f(y)),f(z)))>f(0^{+})$ will lead to a contradiction.
Consequently, $T(T(f(x), f(y)),f(z))),T(c_{k}, f(z))\in(0,f(0^{+}))$, i.e.,
$$(\min\{T(T(f(x), f(y)),f(z)),T(c_{k}, f(z))\},\max\{T(T(f(x), f(y)),f(z)),T(c_{k}, f(z))\})\in (0,f(0^{+})),$$
implying
$$T((\min\{T(f(x), f(y)),c_{k}\},\max\{T(f(x), f(y)),c_{k}\}),f(z))\in (0,f(0^{+})). $$
This follows that
$$T((\min\{T(f(x), f(y)),c_{k}\},\max\{T(f(x), f(y)),c_{k}\}),M\setminus C)\in (0,f(0^{+}))$$ since $f(z)$ is an arbitrary element of $M\setminus C$ satisfying $f(z) > \max\{f(x),f(y)\}$.
Thus
$$T(\cup_{k\in K_1}H_{k}, M\setminus C)\subseteq (0,f(0^{+})).$$
Therefore,
$$T(\cup_{k\in K_1}H_{k}, M\setminus \{0\})\subseteq (0,f(0^{+})),$$
i.e.,
$T(\cup_{k\in K_1}H_{k}, M\setminus\{0\})\cap (M\setminus C)=\emptyset$.
\end{proof}

Clearly, if the function $F$ given by Eq.~\eqref{eq:(5)} is a t-norm then $f$ is a strictly increasing function with $f(1)=1$. Moreover, $T(Q,M)=\emptyset\subseteq  [0,f(0^{+})]$. These together with Propositions \ref{prop4.1}, \ref{prop3.3} and \ref{th4.1} imply the following proposition.
\begin{proposition}\label{th4.02}
Let $T:[0,1]^{2}\rightarrow [0,1]$ be a strict t-norm, $(\mathcal{S},C)$ be associated with $M \in \mathcal{A}^*$ and $\emph{Ran}(f)=M$. Then the function $F$ given by Eq.~\eqref{eq:(5)} is a conditionally cancellative t-norm if and only if the following statements hold:
\renewcommand{\labelenumi}{(\roman{enumi})}
\begin{enumerate}
\item  $f$ is a strictly increasing function with $f(1)=1$,
\item $T(\cup_{k\in K}H_{k}, M\setminus\{0\})\cap (M\setminus C)= \emptyset$ \mbox{and}
\item $T(M\setminus C, M)\subseteq M \cup [0,f(0^{+})]$.
\end{enumerate}
\end{proposition}

\begin{definition}\label{d12.3.1}
	\emph{Let $T:[0,1]^{2}\rightarrow [0,1]$ be a strict t-norm and $(\mathcal{S},C)$ be associated with $M\in \mathcal{N}$. Define for all $y\in M$, $k,l\in K_1$,  $$M_{k}^{y}=\{x\in M \mid T(x,y)\in [b_k, d_k]\},$$ 	$$M^{y}=\{x\in M \mid T(x,y)\in M\setminus C\},$$
	\begin{equation*}
		J_{k,l}^{y}=\begin{cases}
			O(T(M_{k}^{y},c_l)\cup T(c_k,M_{l}^{y})), & \ M_{k}^{y}\neq\emptyset, M_{l}^{y}\neq\emptyset,\\
			\emptyset, & \hbox{otherwise,}
		\end{cases}
	\end{equation*}
 $$I_{k}^{y}=O(\{c_k\}\cup T(M_{k}^{y},y)),$$ 
	and $$\mathfrak{L_{1}}(M)=\bigcup_{y\in M}\bigcup_{k\in K_1}T(I_{k}^{y},M^{y}),$$
	$$\mathfrak{L_{2}}(M)=\bigcup_{y\in M}\bigcup_{k,l\in K_1}J_{k,l}^{y},$$ 
	$$\mathfrak{L}(M)=\mathfrak{L_{1}}(M)\cup\mathfrak{L_{2}}(M).$$}
\end{definition}
Then we have the following lemma.
\begin{lemma}\label{le12.3.2}
	Let $T:[0,1]^{2}\rightarrow [0,1]$ be a strict t-norm and $(\mathcal{S},C)$ be associated with $M\in \mathcal{A}$. Let $M_1, M_2\subseteq [0,1]$ be two non-empty sets and  $c\in[0,1]$. Then
	 $T(O(M_1\cup M_2),c)\cap (M\setminus C)=\emptyset$ if and only if there exist $x_{1}\in M_1$ and $x_{2}\in M_2$ such that
	$T((\min\{x_{1},x_{2}\}, \max\{x_{1},x_{2}\}), c)\cap (M\setminus C)=\emptyset.$
\end{lemma}
\begin{proof}
	 The proof is completely analogous to the proof of Lemma 5.2 of \cite{2025chen}.
\end{proof}

\begin{proposition}\label{prop4.3}
	Let $T:[0,1]^{2}\rightarrow [0,1]$ be a strict t-norm and $(\mathcal{S},C)$ be associated with $M\in \mathcal{N}$. If the function $F$ given by Eq.~\eqref{eq:(5)} is a conditionally cancellative t-subnorm, then $\mathfrak{J}(M)\cap (M\setminus C)=\emptyset$.
\end{proposition}
\begin{proof} Suppose that $\mathfrak{L}(M)\cap (M\setminus C)\neq\emptyset$. Then there exist $y\in M$, $k,l\in K_1$ such that $T(I_{k}^{y},M^{y})\cap (M\setminus C) \neq\emptyset$ or $J_{k,l}^{y}\cap (M\setminus C)\neq \emptyset$. We distinguish two cases as follows.
	
	(i) If $T(I_{k}^{y},M^{y})\cap (M\setminus C) \neq\emptyset$, then there exists a $z\in M^{y}$ such that $$T(I_{k}^{y},z)\cap (M\setminus C) \neq\emptyset.$$ Thus by the definition of $I_{k}^{y}$ , we have $T(M_{k}^{y},y)\neq\emptyset$. From Lemma \ref{le12.3.2}, there exist  $u\in \{c_k\}$, $v\in T(M_{k}^{y},y)$ such that
	$$T((\min\{u,v\}, \max\{u,v\}), z)\cap (M\setminus C)\neq\emptyset.$$
	Moreover, there exists $x\in M_{k}^{y}$ such that $T(x,y)=v$ since 
	$v\in T(M_{k}^{y},y)$. Therefore, there exist $u\in \{c_k\}$, $x\in M_{k}^{y}$ such that
		$$T((\min\{c_k,T(x,y)\}, \max\{c_k,T(x,y)\}), z)\cap (M\setminus C)\neq\emptyset.$$
	Because $f\mid_{(\tau,1]}$ is strictly increasing, by Lemma \ref{lem3.2}, 
	\begin{equation}\label{eq110} 
		f^{(-1)}(T(c_k,z))\neq f^{(-1)}(T(T(x,y),z)).
	\end{equation}
	On the other hand, $x\in M_{k}^{y}$ implies $T(x,y)\in [b_k,d_k]$, thus $f\circ f^{(-1)}(T(x,y))=c_k$. Meanwhile, $z\in M^{y}$ implies $T(y,z)\in M\setminus C$, thus $f\circ f^{(-1)}(T(y,z))=T(y,z)$. Let $x=f(m)$, $y=f(n)$ and $z=f(s)$.
	Then \begin{align*} F(F(m,n),s)&=f^{(-1)}(T(f\circ f^{(-1)}(T(f(m),f(n))),f(s)))\\&=f^{(-1)}(T(f\circ f^{(-1)}(T(x,y)),z))\\&=f^{(-1)}(T(c_k,z))
		\end{align*}
	and \begin{align*}F(m,F(n,s))&=f^{(-1)}(f(m), f\circ f^{(-1)}(T(f(n),f(s))))\\&=
	f^{(-1)}(T(x,f\circ f^{(-1)}(T(y,z))))\\&=f^{(-1)}(T(x,T(y,z))).
	\end{align*}
	Therefore, by Eq.(\ref{eq110}), $F(F(m,n),s)\neq F(m,F(n,s))$, which contradicts the fact that $F$ is a t-subnorm.
	
	(ii) If $J_{k,l}^{y}\cap (M\setminus C)\neq \emptyset$, then $J_{k,l}^{y}\neq \emptyset$. Thus by the definition of $J_{k,l}^{y}$, $$T(O(T(M_{k}^{y},c_l)\cup T(c_k,M_{l}^{y})),1)\cap (M\setminus C)\neq \emptyset,$$ hence $T(M_{k}^{y},c_l)\neq\emptyset$ and $T(c_k,M_{l}^{y})\neq\emptyset$. From Lemma \ref{le12.3.2}, there exist $v\in T(M_{k}^{y},c_l)$ and $u\in T(c_k,M_{l}^{y})$ such that $$T((\min\{u,v\}, \max\{u,v\}), 1)\cap (M\setminus C)\neq\emptyset.$$
	Because $v\in T(M_{k}^{y},c_l)$ and $u\in T(c_k,M_{l}^{y})$, there exist $x\in M_{k}^{y}$ and $z\in M_{l}^{y}$ such that $v=T(x,c_l)$ and $u=T(c_k,z)$. Thus
	$$(\min\{T(c_k,z),T(x,c_l)\},\max\{T(c_k,z),T(x,c_l)\})\cap (M\setminus C)\neq\emptyset$$ since $1$ is the identity element of $T$.
	Hence, by Lemma \ref{lem3.2}, $$f^{(-1)}(T(c_k,z))\neq f^{(-1)}(T(x,c_l)).$$
	On the other hand, $x\in M_{k}^{y}$ implies $T(x,y)\in [b_k,d_k]$, thus $f\circ f^{(-1)}(T(x,y))=c_k$. Meanwhile,
	$z\in M_{l}^{y}$ implies $T(y,z)\in [b_l,d_l]$, thus $f\circ f^{(-1)}(T(y,z))=c_l$. Let $x=f(m)$, $y=f(n)$ and $z=f(s)$.
	Therefore, $$ F(F(m,n),s)=f^{(-1)}(T(c_k,z))\neq f^{(-1)}(T(x,c_l))=F(m,F(n,s)),$$ which contradicts the fact that $F$ is a t-subnorm.
\end{proof}

\begin{proposition}\label{th4.3}
	Let $T:[0,1]^{2}\rightarrow [0,1]$ be a strict t-norm, $(\mathcal{S}, C)$  be associated with  $M \in \mathcal{N}$ and $\emph{Ran}(f)=M$. Then the function $F$ given by Eq.~\eqref{eq:(5)} is a conditionally cancellative t-subnorm if and only if the following conditions hold:
	\renewcommand{\labelenumi}{(\roman{enumi})}
	\begin{enumerate}
		\item $\mathfrak{L}(M)\cap (M\setminus C)=\emptyset$,
		\item $T(M\setminus C, M)\subseteq M \cup [0,f(0^{+})]$ \mbox{and}
		\item $T(Q,M)\subseteq [0,f(0^{+})]$.
	\end{enumerate}
\end{proposition}
\begin{proof}
From Propositions \ref{prop4.1} and \ref{prop4.3} the sufficient condition is obviously.

 Conversely, if $\mathfrak{L}(M)\cap (M\setminus C)=\emptyset$, then $T(\cup_{k\in K_1}H_{k}, M\setminus\{0\})\cap (M\setminus C)= \emptyset$. Therefore, from Proposition \ref{prop4.2} the function $F$ given by Eq.\eqref{eq:(5)} is a conditionally cancellative t-subnorm.
\end{proof}

\begin{proposition}\label{th3.2}
Let $T:[0,1]^{2}\rightarrow [0,1]$ be a strict t-norm, $(\mathcal{S},C)$ be associated with $M \in \mathcal{A}$ and $\emph{Ran}(f)=M$. Then the following are equivalent:
\renewcommand{\labelenumi}{(\roman{enumi})}
\begin{enumerate}
\item The function $F$ given by Eq.\eqref{eq:(5)} is a cancellative t-subnorm.
\item $T(M,M)\subseteq M$.
\end{enumerate}
\end{proposition}
\begin{proof}
 From Corollary \ref{cor.3.2}, (ii) implies (i). Below, we show that (i) implies (ii). Indeed, let $F$ be a cancellative t-subnorm. Supposing $T(M, M) \nsubseteq M $, it is clear that there are two elements $x,y \in [0,1]$ such that $T(f(x),f(y)) \notin M$. Thus there exists a $k \in K$ such that 
 \begin{equation}\label{eq0003}
 	T(f(x),f(y))\in [b_{k}, d_{k}]\setminus\{c_{k}\}
 \end{equation} 
where $\{c_{k}\}=[b_{k}, d_{k}] \cap M$. Note that $(M\setminus C)\setminus\{0\}\neq\emptyset$ since $f$ is strictly increasing. This means that there is a $z\in [0,1]$ such that 
\begin{equation}\label{eq0004}
	f(z)\in (M\setminus C)\setminus \{0\}.
\end{equation}
 On the other hand, by Proposition \ref{prop3.2} we have
\begin{equation}\label{eq1001}
	T(M\setminus C, M) \subseteq (M\setminus C)\cup \{0\}.
\end{equation}
Hence, from Eq.\eqref{eq0004} we have $T(f(y),f(z))\in (M\setminus C)\cup \{0\}$. This deduces that $T(f(y),f(z))\in M\setminus C$ since $f(y)>0,f(z)>0$ and $T$ is a strict t-norm. Thus by Eq.~\eqref{eq:(5)} we have $F(x,F(y,z))=f^{(-1)}(T(f(x), T(f(y),f(z))))$. Meanwhile, Eqs.~\eqref{eq:(5)} and \eqref{eq0003} yield $F(F(x,y),z)=f^{(-1)}(T(c_{k}, f(z)))$. Thereby, the associativity of $F$ and $T$ results in
$$f^{(-1)}(T(T(f(x), f(y)),f(z)))=f^{(-1)}(T(c_{k}, f(z))),$$
implying $$M\cap[\min\{T(T(f(x), f(y)),f(z)),T(c_{k}, f(z))\},\max\{T(T(f(x), f(y)),f(z)),T(c_{k}, f(z))\}]$$ contains at most one element.
Note that from Lemma \ref{lem04.1} we know that $f(0)=0$, which together with Eq. \eqref{eq1001} means that $T(c_{k},f(z))\in (M\setminus C)\cup\{0\}\subseteq M$ and $T(T(f(x), f(y)),f(z)) \in (M\setminus C)\cup\{0\}\subseteq M$. Therefore, $T(T(f(x), f(y)),f(z))=T(c_{k},f(z))$. So that $T(f(x),f(y))=c_{k}$ since $T$ is a strict t-norm and $f(z)>0$, contrary to Eq.\eqref{eq0003}.
\end{proof}

From Proposition \ref{prop4.01}, it is clear that if $f$ is a non-decreasing function that satisfies either $f(1^-)\in Q$ or $f(1)\in Q$ then the function $F$ given by Eq.~\eqref{eq:(5)} is a conditionally cancellative t-subnorm if and only if $T(Q,M)\subseteq [0,f(0^{+})]$. Therefore, by Propositions \ref{prop4.01}, \ref{th4.3} and \ref{th3.2} the following two theorems are immediate.
\begin{theorem}\label{th4.4}
		Let $T:[0,1]^{2}\rightarrow [0,1]$ be a strict t-norm, $(\mathcal{S},C)$ be associated with $M \in \mathcal{A}^*$ and $\emph{Ran}(f)=M$. Then the function $F$ given by Eq.~\eqref{eq:(5)} is a conditionally cancellative t-subnorm if and only if one of the following holds:
		\renewcommand{\labelenumi}{(\roman{enumi})}
		\begin{enumerate}
			\item  $f$  is a non-increasing function with $f(x)=0$ for all $x\in(0,1]$.
			\item  $f$  is a non-decreasing function satisfying
			that \renewcommand{\labelenumi}{(\roman{enumi})}
			\begin{enumerate}
				\item $\mathfrak{L}(M)\cap (M\setminus C)=\emptyset$, 
				\item $T(M\setminus C, M)\subseteq M \cup [0,f(0^{+})]$ \mbox{and}
				\item $T(Q,M)\subseteq [0,f(0^{+})]$.
			\end{enumerate}
		\end{enumerate}
\end{theorem}

\begin{theorem}\label{th3.02}
Let $T:[0,1]^{2}\rightarrow [0,1]$ be a strict t-norm and $f:[0,1]\rightarrow [0,1]$ be a monotone function with $\emph{Ran}(f)=M$. Then the function $F$ given by Eq.~\eqref{eq:(5)} is a cancellative t-subnorm if and only if the following two statements hold:
\renewcommand{\labelenumi}{(\roman{enumi})}
\begin{enumerate}
\item $f$ is a strictly increasing function \mbox{and}
\item $T(M,M)\subseteq M$.
\end{enumerate}
\end{theorem}
\begin{example}\label{exap3.1}
\renewcommand{\labelenumi}{(\roman{enumi})}
\emph{\begin{enumerate}
\item Let the function $f:[0,1]\rightarrow [0,1]$ be defined by
		\begin{equation*}
			f(x)=\begin{cases}
				\frac{1}{2} & \hbox{if }\ x\in[0,\frac{1}{2}],\\
				x & \hbox{if }\ x\in(\frac{1}{2},1]
			\end{cases}
		\end{equation*}
		and  $T:[0,1]^{2}\rightarrow [0,1]$ be defined by
			$T(x,y)=xy$.
		Then by Eq.~\eqref{eq:(5)}
		\begin{equation*}
			F(x,y)=\begin{cases}
				0 & \hbox{if }\ 0<xy\leq \frac{1}{2},\\
				xy & \hbox{otherwise}.
			\end{cases}
		\end{equation*}
		It is easy to see that $T$ is a strict t-norm. One can check that $\mbox{Ran}(f)= M=[\frac{1}{2},1]$ and $T(f(x),f(y))\in M \cup [0,f(0^{+})]$ for all $(x,y)\in[0,1]^{2}$ and $T(Q,M)\subseteq [0,f(0^{+})]$. Therefore, from Theorem \ref{th4.4} $F$ is a conditionally cancellative t-subnorm.
\item Let the function $f:[0,1]\rightarrow [0,1]$ be defined by
		\begin{equation*}
			f(x)=\begin{cases}
				\frac{x}{2} & \hbox{if }\ x\in[0,1),\\
				1 & \hbox{if }\ x=1
			\end{cases}
		\end{equation*}
		and  $T:[0,1]^{2}\rightarrow [0,1]$ be defined by
		\begin{equation*}
			T(x,y)=\frac{xy}{2-(x+y-xy)}.
		\end{equation*}
		Then by Eq.~\eqref{eq:(5)}
		\begin{equation*}
			F(x,y)=\begin{cases}
				\frac{xy}{8+xy-2(x+y)} & \hbox{if }\ (x,y)\in[0,1)^{2},\\
				\min\{x,y\} & \hbox{otherwise}.
			\end{cases}
		\end{equation*}
		It is easy to see that $f$ is strictly increasing with $\mbox{Ran}(f)=[0,\frac{1}{2})\cup\{1\}$. One can check that $T$ is a strict t-norm and $T(f(x),f(y))\in \mbox{Ran}(f)$ for all $(x,y)\in[0,1]^{2}$. Therefore, from Theorem \ref{th3.02} $F$ is a cancellative t-subnorm.
\end{enumerate}}
\end{example}

From Theorem \ref{th3.02} we have the following remark.
\begin{remark}\label{re3.1}
\renewcommand{\labelenumi}{(\roman{enumi})}
\emph{\begin{enumerate}
\item Let $T:[0,1]^{2}\rightarrow [0,1]$ be a strict t-norm and $f:[0,1]\rightarrow [0,1]$ be a monotone function with $\mbox{Ran}(f)=M$. Then the function $F:[0.1]^{2}\rightarrow [0,1]$, given by
\begin{equation*}
F(x,y)=\begin{cases}
f^{(-1)}T(f(x),f(y)) & \hbox{if }\ (x,y)\in[0,1)^{2},\\
\min\{x,y\} &  \hbox{otherwise},\\
\end{cases}
\end{equation*}
is a cancellative t-norm if and only if $f$ is strictly increasing and $T(M\setminus f(1),M\setminus f(1))\subseteq M\setminus f(1)$.
\item Let $T:[0,1]^{2}\rightarrow [0,1]$ be a strict t-norm and $f:[0,1]\rightarrow [0,1]$ be a monotone function with $\mbox{Ran}(f)=M$. Then the function $F$ given by Eq.\eqref{eq:(5)} is a cancellative t-norm if and only if $f$ is a strictly increasing function with $f(1)=1$ and $T(M,M)\subseteq M$. In particular, let $T(x,y)=x\cdot y$. Then $F(x,y)=f^{(-1)}(f(x)\cdot f(y))$ is a cancellative t-norm if and only if $f$ is strictly increasing with $f(1)=1$ and
    $f(x)\cdot f(y)\in \mbox{Ran}(f)$ for all $x,y \in[0,1]$.
\end{enumerate}}
\end{remark}

As a conclusion of this section, we point out that if $f:[0,1]\rightarrow [0,1]$ is a continuous, strictly increasing function and $T$ is a strict t-norm then from Theorems \ref{th2.1} and \ref{th2.2}, the function $F$ given by Eq.\eqref{eq:(5)} is a strict t-subnorm. However,
the converse is not true generally. For instance, let the function
$T:[0, 1]^{2} \rightarrow [0, 1]$ be given by \begin{equation*}
T(x,y)=\begin{cases}
\frac{xy}{x+y} & \hbox{if }\ (x,y)\in[0,1)^{2},\\
\min\{x,y\} &  \hbox{otherwise}
\end{cases}
\end{equation*}
and the function $f:[0, 1] \rightarrow [0, 1]$ be defined by $f(x)=\frac{x}{2}$. Then by Eq.~\eqref{eq:(5)}
$$F(x,y)=\frac{xy}{x+y} \hbox{ for all }\ (x,y)\in[0,1]^{2}.$$
It is easily checked that $F$ is a strict t-subnorm and $T$ is not a continuous t-norm.

Fortunately, we have the following theorem.
 \begin{theorem}\label{th3.4}
For a function $F: [0, 1]^{2} \rightarrow [0, 1]$ the following are equivalent:
\renewcommand{\labelenumi}{(\roman{enumi})}
\begin{enumerate}
\item $F$ is a strict proper t-subnorm.
\item There exist a continuous,
strictly increasing function $f: [0, 1] \rightarrow [0, 1]$ with $f(0)=0$ and $f(1)<1$, and a
strictly monotone t-norm $T: [0, 1]^{2} \rightarrow [0, 1]$ such that
Eq.~\eqref{eq:(5)} holds.
\end{enumerate}
\end{theorem}
\begin{proof}
From Corollary \ref{cor.3.3}, (ii) implies (i). We prove that (i) implies (ii). Let $F$ be a strict t-subnorm. Then from Theorem \ref {th2.2} there exists a continuous,
strictly decreasing function $g:[0, 1] \rightarrow [0, \infty]$ with $g(0)=\infty$ and $g(1)>0$ such that
$F(x,y) = g^{(-1)}(g(x)+g(y))$ for all $(x , y) \in [0, 1]^{2}$.
For any $\lambda\in(0,1)$, define a function $t:[0, 1] \rightarrow [0, \infty]$ by
\begin{equation*}
t(x)=\begin{cases}
g(\frac{x}{\lambda}) & \hbox{if }\ x\in[0,1),\\
0 &  \hbox{if }\ x=1.\\
\end{cases}
\end{equation*}
Then the function $T:[0, 1]^{2} \rightarrow [0, 1]$, given by
\begin{equation*}
T(x,y)=\begin{cases}
t^{-1}(t(x)+t(y)) & \hbox{if }\ (x,y)\in[0,1)^{2},\\
\min\{x,y\} &  \hbox{otherwise},
\end{cases}
\end{equation*}
is a strictly monotone t-norm.
Define a function $f:[0, 1] \rightarrow [0, 1]$ by $f(x)=\lambda x$. It is easily checked that $f$ is a continuous,
strictly increasing function with $f(0)=0$ and $f(1)<1$, and $g=t\circ f$ where $(t\circ f)(x)=t(f(x))$ for all $x\in [0,1]$. Thus
\begin{align*}
F(x,y) &= g^{(-1)}(g(x)+g(y)) \\
&=f^{(-1)}\circ t^{-1} (t\circ f(x)+t\circ f(y)) \\
&=f^{(-1)}(T(f(x) ,f(y)))
\end{align*}
 for all $(x,y) \in[0, 1]^{2}$.
\end{proof}
\section{Conclusions}
This artcile presented a characterizations of monotone function $f$ such that the function $F: [0,1]^2\rightarrow [0,1]$ given by Eq.\eqref{eq:(5)} is a conditionally cancellative t-subnorm (resp. a cancellative t-subnorm) (Theorem \ref{th4.4} (resp. Theorem \ref{th3.02})). We also characterized that a function $F$ given by Eq.\eqref{eq:(5)} is a strict proper t-subnorm (Theorem \ref{th3.4}). It is easy to see that we can give a characterizations of monotone function $f$ such that the function $F: [0,1]^2\rightarrow [0,1]$ given by Eq.\eqref{eq:(5)} is a conditionally cancellative t-supconorm (resp. a cancellative t-supconorm). For instance, a function $F$ given by Eq.\eqref{eq:(5)} is a conditionally cancellative t-supconorm if and only if $f$ is a non-increasing function, $\mathfrak{L}(M)\cap (M\setminus C)=\emptyset$, $T(M\setminus C, M)\subseteq M \cup [0,f(1^{-})]$ and $T(Q,M)\subseteq [0,f(1^{-})]$. 

In particular, from our results we can answer the open problem posed by Mesiarov\'{a} \cite{AM2004}. Indeed, because $f\mid_{[\tau,1]}$ is strictly increasing where $\tau \in [0,1]$ with $f(\tau^{-})=\upsilon$, it follows from Lemma \ref{lemma2.1}, Proposition \ref{prop2.1} and Theorem \ref{th4.4} that we have the following theorem.
\begin{theorem}\label{th5.1}
	Let $T:[0,1]^{2}\rightarrow [0,1]$ be a strict t-norm, $(\mathcal{S},C)$ be associated with $M \in \mathcal{A}^*$ and $\emph{Ran}(f)=M$. Then the function $F$ given by Eq.~\eqref{eq:(5)} is a continuous conditionally cancellative proper t-subnorm if and only if one of the following holds:
	\renewcommand{\labelenumi}{(\roman{enumi})}
	\begin{enumerate}
		\item $f$  is a non-increasing function with $f(x)=0$ for all $x\in(0,1]$.
		\item $f$ is a non-decreasing function satisfying that 
		\renewcommand{\labelenumi}{(\roman{enumi})}
		\begin{enumerate}
			\item $\emph{Ran}(f)\cap[T(f(x^-),f(y^-)),T(f(x^+),f(y^+))]$ is at most a one-element set for all $x,y\in(\tau,1]$ and, $f(1)<1$ whenever $f$ is strictly increasing,
			\item $ \mathfrak{L}(M)\cap (M\setminus C)=\emptyset$, 
			\item $T(M\setminus C, M)\subseteq M \cup [0,f(0^{+})]$ \mbox{and}
			\item  $T(Q,M)\subseteq [0,f(0^{+})]$.
		\end{enumerate}
	\end{enumerate}
\end{theorem}

On the one hand, because of the complete duality between additive and multiplicative generators, the open problem posed by Mesiarov\'{a} \cite{AM2004} can be restated as: Characterize all multiplicative generators yielding continuous Archimedean proper t-subnorms? On the other hand, it is easy to verify that 
a continuous t-subnorm is Archimedean if and only if it satisfies the conditional cancellation law. Thus, the above open problem can further be rephrased as: Characterize all multiplicative generators yielding continuous conditionally cancellative proper t-subnorms? Therefore, by Theorem \ref{th5.1} we characterize all multiplicative generators yielding continuous conditionally cancellative proper t-subnorms, answering the open problem.


\begin{thebibliography}{99}
\bibitem{Abel}
N.H. Abel, Untersuchung der Function zweier unabhangig veranderlichen Grossen $x$ und $y$ wie $f(x,y)$, welche die Eigenschaft haben, dass $f(z,f(x,y))$ eine symmetrische Function von $x$, $y$ und $z$ ist, J. Reine Angew. Math. 1 (1826) 11-15.

\bibitem{CA2006}
 C. Alsina, M.J. Frank, B. Schweizer, Associative Functions: Triangular Norms and Copulas, World Scientific, Singapore, 2006.

\bibitem{2025chen}
M. Chen, Y.M. Zhang, X.P. Wang, The characterization of monotone functions that generate associative functions, Fuzzy Sets Syst. 500 (2025) 109201.
\bibitem{Jenei} S. Jenei, A note on the ordinal sum theorem and its consequence for the construction of triangular norms, Fuzzy Sets Syst. 126 (2002) 199–205.
\bibitem{EP2000}
 E.P. Klement, R. Mesiar, E. Pap, Triangular Norms, Kluwer Academic Publishers, Dordrecht, 2000.
\bibitem{Klement} E.P. Klement, R. Mesiar, E. Pap, Triangular norms as ordinal sums of semigroups in the sense of A.H. CliKord, Semigroup Forum 65 (2002) 71–82.
\bibitem{AM2002}
A. Mesiarov\'{a}, Triangular norms and their diagonal functions, Faculty of Mathematics Physics and Informatics UK, Bratislava, 2002 Master thesis.

\bibitem{AM2004}
A. Mesiarov\'{a}, Continuous triangular subnorms, Fuzzy Sets Syst. 142 (2004) 75-83.

\bibitem{AM2016}
A. Mesiarov\'{a}-Zem\'{a}nkov\'{a}, Continuous additive generators of continuous, conditionally cancellative triangular subnorms, Inf. Sci. 339 (2016) 53-63.


\bibitem{PV1999}P. Vicen\'{\i}k, A note to a construction of t-norms based on pseudo-inverses of monotone functions, Fuzzy Sets Syst. 104 (1999) 15-18.

\bibitem{PV2005}
P. Vicen\'{\i}k, Additive generators of associative functions, Fuzzy Sets Syst. 153 (2005) 137-160.

\bibitem{Wang2005}
	X.P. Wang, Y.M. Zhang, Some algebraic and analytical properties of a class of two-place functions, Fuzzy Sets Syst. 500 (2025) 109196.
\bibitem{YM2024}
Y.M. Zhang, X.P. Wang, Characterizations of monotone right continuous functions which generate associative functions, Fuzzy Sets Syst. 477 (2024) 108799.
\end{thebibliography}
\end{document}